\documentclass[11pt,a4paper]{article}
\usepackage{amsfonts}
\usepackage{latexsym}
\usepackage{cite}
\usepackage{amsmath,amsfonts,latexsym,amssymb}
\usepackage[mathscr]{eucal}
\usepackage{cases,color}
\usepackage{amsthm}

\usepackage[bf,small]{caption2}
\usepackage{float}
\usepackage{graphicx}
\usepackage{amsmath}
\usepackage{amssymb}
\usepackage[all]{xy}

\newtheorem{theorem}{theorem}[section]

\newtheorem{lem}[theorem]{Lemma}
\newtheorem{nota}[theorem]{Notation}
\newtheorem{prob}[theorem]{Problem}

\newtheorem{rmk}[theorem]{Remark}
\newtheorem{thm}[theorem]{Theorem}

\begin{document}

\title{\vspace{-2cm}\textbf{Kauffman bracket skein algebra of the 4-holed disk}}
\author{\Large Haimiao Chen}
\date{}
\maketitle

\begin{abstract}
  We give a monomial basis for the Kauffman bracket skein algebra of the $4$-holed disk, and find a presentation.
  This is based on an insight into the ${\rm SL}(2,\mathbb{C})$-character variety of the rank $4$ free group.

\end{abstract}

\section{Introduction}

Let $R$ be a commutative ring with identity and a fixed invertible element $q^{\frac{1}{2}}$.
For an oriented $3$-manifold $M$, its {\it Kauffman bracket skein module} over $R$, denoted by $\mathcal{S}(M;R)$, is defined as the quotient of the free $R$-module generated by isotopy classes of (possibly empty) framed links embedded in $M$ by the submodule generated by the following {\it skein relations}:
\begin{figure}[h]
  \centering
  \includegraphics[width=9cm]{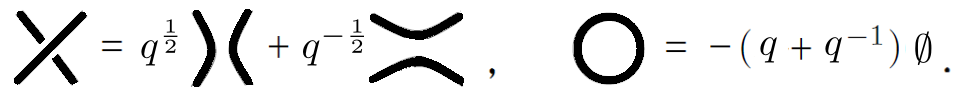}\\
  %\caption{}\label{fig:local}
\end{figure}
\\
As a convention, $R$ is identified with $R\emptyset\subset\mathcal{S}(M;R)$ via $\mu\mapsto \mu\emptyset$.

When $M=\Sigma\times[0,1]$ for an oriented surface $\Sigma$, denote $\mathcal{S}(M;R)$ as $\mathcal{S}(\Sigma;R)$ and call it the {\it skein algebra} of $\Sigma$.
Given links $\mathbf{l}_1,\mathbf{l}_2\subset\Sigma\times[0,1]$, the product $\mathbf{l}_1\mathbf{l}_2$ is defined by stacking $\mathbf{l}_1$ over $\mathbf{l}_2$ in the $[0,1]$ direction. Adopt the convention that each framed link is presented as a link equipped with the blackboard framing, i.e. each framing vector points vertically upward.

Let $\Sigma_{g,k}$ denote a $k$-holed genus $g$ orientable surface.
The following was raised as \cite[Problem 1.92 (J)]{Ki97} and \cite[Problem 4.5]{Oh02}:
\begin{prob}[Bullock and Przytycki]\label{prob}
Find the structure of $\mathcal{S}(\Sigma_{g,k};\mathbb{Z}[q^{\pm\frac{1}{2}}])$.
\end{prob}

The structure of $\mathcal{S}(\Sigma_{g,k};\mathbb{Z}[q^{\pm\frac{1}{2}}])$ for $g=0,k\le 4$ and $g=1,k\le 2$ was known to Bullock and Przytycki \cite{BP00} early in 2000. A finite set of generators had been given by Bullock \cite{Bu99} in 1999. Till now it remains a difficult problem to find all relations for general $g$ and $k$.

Recently, Cooke and Lacabanne \cite{CL22} gave a presentation for $\mathcal{S}(\Sigma_{0,5};\mathbb{C}(q^{\frac{1}{4}}))$, and the author \cite{Ch24} found a presentation for $\mathcal{S}(\Sigma_{0,n+1};R)$ for all $n\ge 4$ and all $R$ containing the inverse of $q+q^{-1}$.

Let $\mathcal{S}_n=\mathcal{S}(\Sigma_{0,n+1};\mathbb{Z}[q^{\pm\frac{1}{2}}])$. In this paper we focus on $\mathcal{S}_4$.

\begin{figure}[H]
  \centering
  \includegraphics[width=8.5cm]{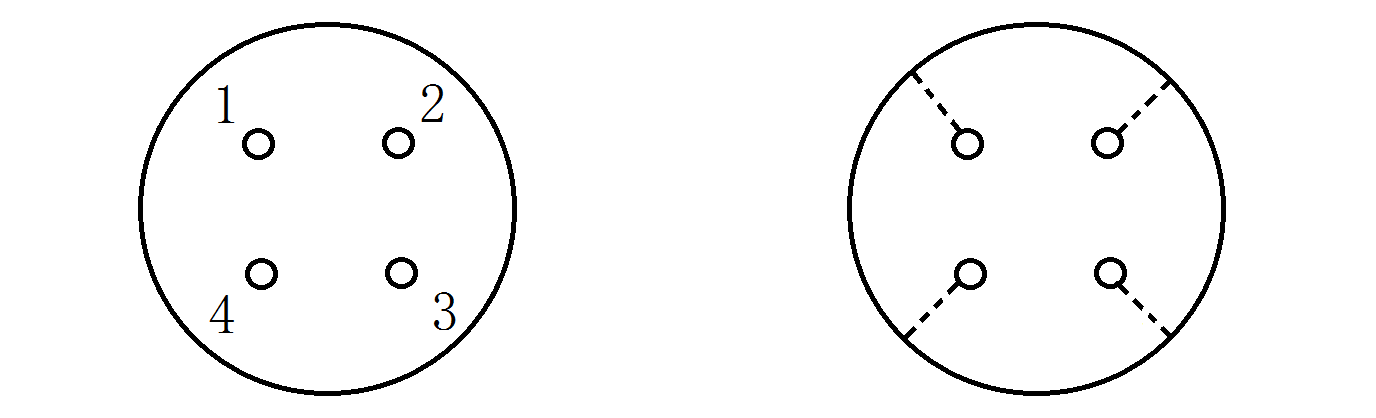}\\
  \caption{Left: the $4$-holed disk $\Sigma_{0,5}$. Right: $\mathbf{z}_1,\ldots,\mathbf{z}_4$.}\label{fig:Sigma}
\end{figure}

Display $\Sigma_{0,5}$ as in Figure \ref{fig:Sigma}. We prefer to call $\Sigma_{0,5}$ the $4$-holed disk rather than the $5$-holed sphere, because the outer boundary will play a distinct role. Let $\mathbf{z}_i$ be an arc connecting the $i$-th hole to the outer boundary.

\begin{figure}[H]
  \centering
  \includegraphics[width=12.5cm]{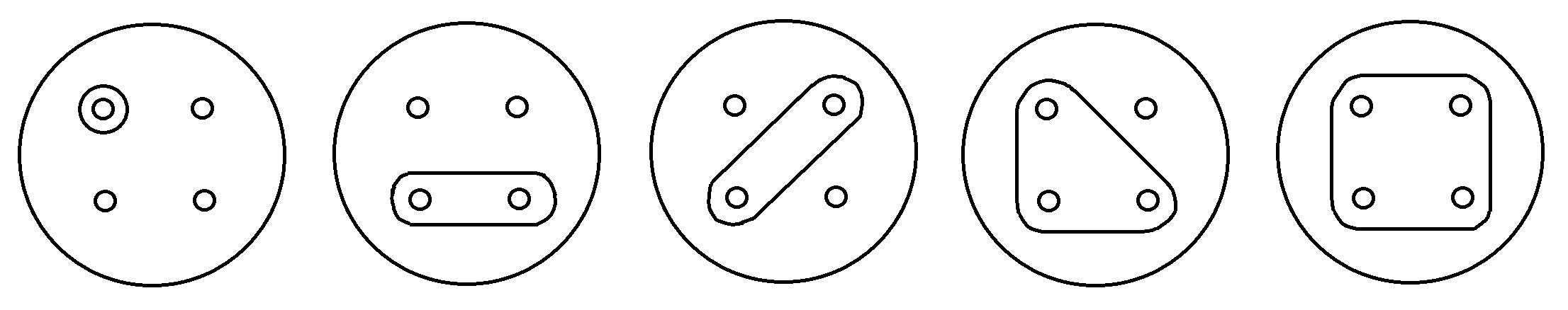}\\
  \caption{From left to right: $t_1$, $t_{34}$ $t_{24}$, $t_{134}$, $t_0=t_{1234}$.}\label{fig:notation-1}
\end{figure}

For $1\le i_1<\cdots<i_r\le 4$, let $t_{i_1\cdots i_r}$ be a simple curve intersecting $\mathbf{z}_j$ once exactly for $j\in\{i_1,\ldots,i_r\}$; see Figure \ref{fig:notation-1} for some instances.
Let
$$\mathcal{G}=\{t_{i_1\cdots i_r}\colon 1\le i_1<\cdots<i_r\le 4\}.$$
Denote $t_{1234}$ by $t_0$.
Note that $t_1,t_2,t_3,t_4,t_0$ are central in $\mathcal{S}_4$.
Call a product of elements of $\mathcal{G}$ a {\it monomial}, and call a $R$-linear combination of monomials a {\it polynomial}.

Let the cyclic permutation $\sigma=(1234)$ act on $\mathcal{G}$ by
$\sigma(t_{i_1\cdots i_r})=t_{\sigma(i_1)\cdots\sigma(i_r)}$, and extend the action to polynomials $f$.
If $f=0$ is a relation, then so is $\sigma^k(f)=0$ for $1\le k\le 3$; call it the relation obtained from rotating $f=0$ by $-k\pi/2$.

\begin{nota}
\rm Denote $q^{-1}$ by $\overline{q}$, denote $q^{-\frac{1}{2}}$ by $\overline{q}^{\frac{1}{2}}$, and so on. Let $\alpha=q+\overline{q}$.

For a finite set $A$, let $\#A$ denote its cardinality.

Given a subset $X$ of a left module $\mathcal{M}$ over a ring $\mathcal{R}$, let $\mathcal{R}\langle X\rangle$ denote the submodule of $\mathcal{M}$ generated by $X$.
\end{nota}

In this paper, we establish the following two theorems.
\begin{thm}\label{thm:main-1}
As a $\mathbb{Z}[q^{\pm\frac{1}{2}}][t_1,t_2,t_3,t_4]$-module, $\mathcal{S}_4$ is freely generated by
$$\big\{t_{12}^{j_1}t_{23}^{j_2}t_{34}^{j_3}t_{14}^{j_4}a\colon j_1,\ldots,j_4\ge 0,\ a\in\mathcal{A}\big\},$$
where $\mathcal{A}$ consists of $1,t_0,t_{123},t_{124},t_{134},t_{234}$ and
\begin{align*}
t_{13}^k,\ \ t_{13}^kt_0,\ \ t_{13}^kt_{123},\ \  t_{13}^kt_{134}, \ \ t_{24}^k, \ \ t_{24}^kt_0, \ \ t_{24}^kt_{124}, \ \ t_{24}^kt_{234} \ \ \ \text{for\ \ }k\ge 1.
\end{align*}
\end{thm}

\begin{thm}\label{thm:main-2}
The skein algebra $\mathcal{S}_4$ %over any commutative ring containing $q^{\pm\frac{1}{2}}$
is generated by $\mathcal{G}$, and the ideal of defining relations is generated by $\mathcal{H}$ consisting of the following three families of relations.

(a) Commuting relations: $t_1$, $t_2$, $t_3$, $t_4$, $t_0$ are central;
\begin{align*}
t_{34}t_{12}=t_{12}t_{34},  \quad  t_{123}t_{12}=t_{12}t_{123}, \quad
t_{124}t_{12}=t_{12}t_{124}, \quad  t_{123}t_{13}=t_{13}t_{123},
\end{align*}
and the ones obtained by rotations.

(b) Commutator relations:
\begin{align}
t_{23}t_{12}&=\overline{q}^2t_{12}t_{23}+(q-\overline{q}^3)t_{13}+(1-\overline{q}^2)(t_1t_3+t_2t_{123}), \label{eq:commutator-23-12}  \\
t_{13}t_{12}&=q^2t_{12}t_{13}+(\overline{q}-q^3)t_{23}+(1-q^2)(t_2t_3+t_1t_{123}),  \label{eq:commutator-13-12} \\
t_{24}t_{12}&=\overline{q}^2t_{12}t_{24}+(q-\overline{q}^3)t_{14}+(1-\overline{q}^2)(t_1t_4+t_2t_{124}), \label{eq:commutator-24-12}  \\
t_{234}t_{14}&=q^2t_{14}t_{234}+(\overline{q}-q^3)t_{123}+(1-q^2)(t_4t_0+t_1t_{23}), \label{eq:commutator-234-14} \\
t_{124}t_{34}&=\overline{q}^2t_{34}t_{124}+(q-\overline{q}^3)t_{123}+(1-\overline{q}^2)(t_4t_0+t_3t_{12}),  \label{eq:commutator-124-34}
\end{align}
and the ones obtained by rotations.

(c) Reduction relations:
\begin{align}
t_{13}t_{24}&=\alpha t_0+t_1t_{234}+t_2t_{134}+t_3t_{124}+t_4t_{123}+q^2t_{12}t_{34}+\overline{q}^2t_{14}t_{23} \nonumber  \\
&\ \ \ \ +qt_3t_4t_{12}+\overline{q}t_1t_4t_{23}+qt_1t_2t_{34}+\overline{q}t_2t_3t_{14}+t_1t_2t_3t_4,   \label{eq:reduction-13-24}  \\
t_{24}t_{134}&=t_{234}t_{14}+t_{124}t_{34}-t_4t_0-\alpha t_{123}-t_1t_{23}-t_3t_{12}  \nonumber  \\
&\ \ \ \ +t_2(\overline{q}t_{34}t_{14}-\overline{q}^2t_{13}-\overline{q}t_4t_{134})-\overline{q}t_1t_2t_3,  \label{eq:reduction-24-134}   \\
t_{134}t_{24}&=t_{14}t_{234}+t_{34}t_{124}-t_4t_0-\alpha t_{123}-t_1t_{23}-t_3t_{12}  \nonumber  \\
&\ \ \ \ +t_2(qt_{14}t_{34}-q^2t_{13}-qt_4t_{134})-qt_1t_2t_3,   \label{eq:reduction-134-24}  \\
t_{123}^2&=\overline{q}t_{12}t_{23}t_{13}-(t_1t_2t_3+qt_1t_{23}+\overline{q}t_2t_{13}+\overline{q}t_3t_{12})t_{123}-t_1^2-t_2^2-t_3^2 \nonumber  \\
&\ \ \ \ +\alpha^2-qt_2t_3t_{23}-\overline{q}t_1t_3t_{13}-\overline{q}t_1t_2t_{12}-q^2t_{23}^2-\overline{q}^2t_{13}^2-\overline{q}^2t_{12}^2, \label{eq:reduction-123-123}  \\
t_{123}t_{234}&=(t_{23}+\overline{q}t_2t_3)t_0+\overline{q}t_{12}t_{23}t_{34}-\overline{q}t_3t_{12}t_{234}-\overline{q}t_2t_{123}t_{34}
+\overline{q}^2t_2t_{124}   \nonumber  \\
&\ \ \ \ +\overline{q}^2t_3t_{134}-\overline{q}^2t_{12}t_{24}-\overline{q}^2t_{13}t_{34}+\overline{q}^2(\alpha t_{14}+t_1t_4), \label{eq:reduction-123-234}    \\
t_{234}t_{123}&=(t_{23}+qt_2t_3)t_0+qt_{34}t_{23}t_{12}-qt_3t_{234}t_{12}-qt_2t_{34}t_{123}+q^2t_2t_{124}  \nonumber  \\
&\ \ \ \ +q^2t_3t_{134}-q^2t_{24}t_{12}-q^2t_{34}t_{13}+q^2(\alpha t_{14}+t_1t_4),   \label{eq:reduction-234-123}  \\
t_{123}t_{134}&=t_{13}t_0+t_{12}t_{14}+t_{23}t_{34}-t_1t_{123}-t_3t_{234}-\alpha t_{24}-t_2t_4,  \label{eq:reduction-123-134} \\
t_0t_{234}&=\overline{q}(t_{23}t_{34}-\overline{q}t_{24}-t_2t_4)t_{124}-(qt_2t_{34}+\overline{q}t_3t_{24}+\overline{q}t_4t_{23}+t_2t_3t_4)t_0 \nonumber  \\
&\ \ \ \ -q^2t_{34}t_{134}-\overline{q}^2t_{23}t_{123}-\overline{q}t_3t_{23}t_{12}-qt_3t_{14}t_{34}  \nonumber \\
&\ \ \ \ +qt_1t_3^2+q^2t_3t_{13}-t_2t_{12}-t_4t_{14}-\alpha t_1,  \label{eq:reduction-0-234}  \\
t_0^2&=\overline{q}t_{12}t_{234}t_{134}-(t_1t_2t_{34}+qt_1t_{234}+\overline{q}t_2t_{134}+\overline{q}t_{12}t_{34})t_0  \nonumber \\
&\ \ \ \ -t_1^2-t_2^2-t_{34}^2+\alpha^2-qt_2t_{34}t_{234}-\overline{q}t_1t_{34}t_{134}-\overline{q}t_1t_2t_{12}  \nonumber  \\
&\ \ \ \ -q^2t_{234}^2-\overline{q}^2t_{134}^2-\overline{q}^2t_{12}^2,  \label{eq:reduction-0-0}
\end{align}
and the ones obtained by rotations.
\end{thm}

Thus, we have an answer to Problem \ref{prob} for $g=0,k=5$.

That $\mathcal{G}$ generates $\mathcal{S}_4$ is a special case of \cite[Theorem 1]{Bu99}.
We will verify the relations in $\mathcal{H}$, and then prove Theorem \ref{thm:main-1} by showing (i) elements of
$$\mathcal{C}:=\big\{t_1^{i_1}t_2^{i_2}t_3^{i_3}t_4^{i_4}t_{12}^{i_5}t_{23}^{i_6}t_{34}^{i_7}t_{14}^{i_8}a\colon i_1,\ldots,i_8\ge 0, \ a\in\mathcal{A}\big\}$$
are $\mathbb{Z}[q^{\pm\frac{1}{2}}]$-linearly independent, and (ii) each monomial can be transformed into a $\mathbb{Z}[q^{\pm\frac{1}{2}}]$-linear combination of elements of $\mathcal{C}$ via relations in $\mathcal{H}$. This also proves Theorem \ref{thm:main-2}.

The proof of (i) relies on a new understanding about the structure of the ${\rm SL}(2,\mathbb{C})$-character variety of $F_4$, the free group of rank $4$. We present Theorem \ref{thm:classical-0} and Theorem \ref{thm:classical} which have their own interests.

The presentation for $\mathcal{S}_4$ given in Theorem \ref{thm:main-2} is essentially the same as the one in \cite[Theorem 9.3]{CL22}.
There, the coefficient ring is $\mathbb{C}(q^{\frac{1}{4}})$, and the proof applied sophisticated tools such as Hilbert series, the diamond lemma, and term rewriting system, with the aid of a computer. In contrast, in this paper, the coefficient ring is the original $\mathbb{Z}[q^{\pm\frac{1}{2}}]$ in Problem \ref{prob}, and the proof is explicit, constructive, and relatively elementary.

We believe that our method can be easily extended to $\Sigma_{1,3}$ and $\Sigma_{2,1}$, whose fundamental groups are isomorphic to $F_4$.
Then it will cost moderate efforts to derive a presentation for $\mathcal{S}(\Sigma_{2,0};\mathbb{Z}[q^{\pm\frac{1}{2}}])$, as $\Sigma_{2,0}\times[0,1]$ results from $\Sigma_{2,1}\times[0,1]$ by attaching a $2$-handle.

Nowadays, skein modules have become central objects in quantum topology.
An important reason is that %in light of (\ref{eq:classicalization}),
skein module is regarded as a {\it quantization} of character variety \cite{Bu97}.
Computational results were seen in \cite{AF22,DW21,Le06,LT14,LSW24,Ma10} and the references therein.
To compute skein modules for more $3$-manifolds, the first step is to determine that for a handlebody, which is diffeomorphic to $\Sigma_{0,n+1}\times[0,1]$ if the genus is $n$.
Although the skein algebra $\mathcal{S}_n$ is known to admit a basis given by multicurves in $\Sigma_{0,n+1}$, in practice it is often not convenient to work with multicurves.

The monomial basis given by Theorem \ref{thm:main-1} will be advantageous in algebraic manipulations.
We expect it to be applicable in determining skein modules for $3$-manifolds resulted from attaching $2$-handles to genus $4$ handlebodies.
As one instance, in \cite{Ch25} we apply Theorem \ref{thm:main-1} to compute the skein module of the $(3,3,3,3)$-pretzel link exterior.

\section{Verifying the relations in $\mathcal{H}$}

\begin{figure}[H]
  \centering
  \includegraphics[width=9.5cm]{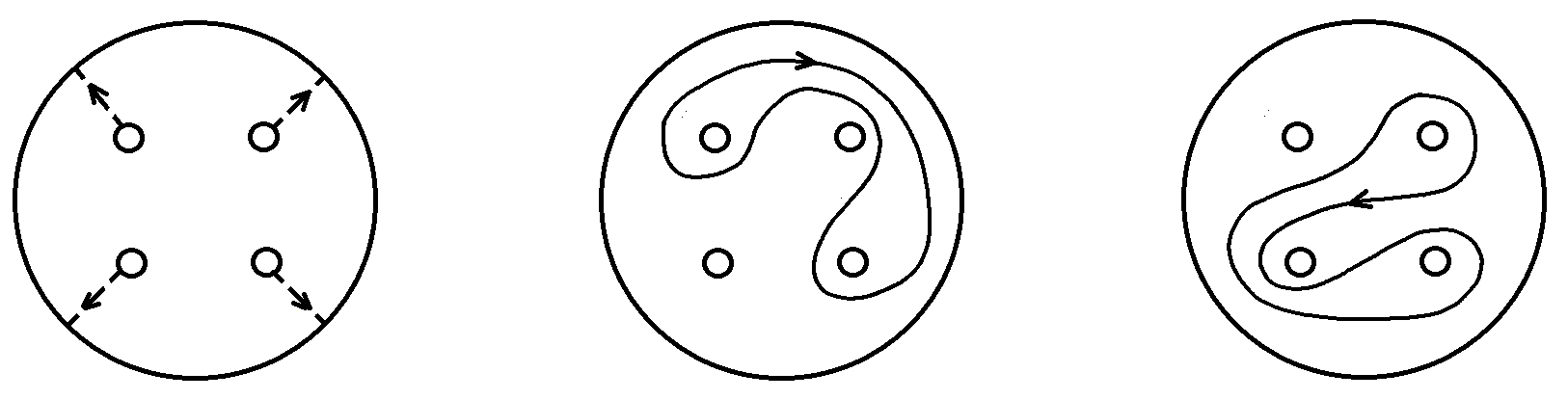}\\
  \caption{Left: $\mathbf{z}_i$ is oriented outward. \ Middle: $t_{123\overline{2}}=t_{3\overline{2}12}$.\ Right: $t_{2\overline{4}34}=t_{342\overline{4}}$.}\label{fig:notation-2}
\end{figure}

Orient each $\mathbf{z}_i$ outward.
Let $\mathbf{s}\subset\Sigma_{0,5}$ be a simple curve.
Starting at a point $\mathsf{x}\in \mathbf{s}$, walk along $\mathbf{s}$ in any direction, and record $i^\ast=i$ (resp. $i^\ast=\overline{i}$) whenever passing through $\mathbf{z}_i$ from left to right (resp. from right to left); the meanings of ``left" and ``right" are self-evident.
Denote $\mathbf{s}$ as $t_{i_1^\ast\cdots i_r^\ast}$ if when back to $\mathsf{x}$, the recordings are $i_1^\ast,\ldots,i_r^\ast$. Depending on the choices of $\mathsf{x}$ and the direction, there are several ways of denoting $\mathbf{s}$. See Figure \ref{fig:notation-2} for examples.

\begin{figure}[H]
  \centering
  \includegraphics[width=11.5cm]{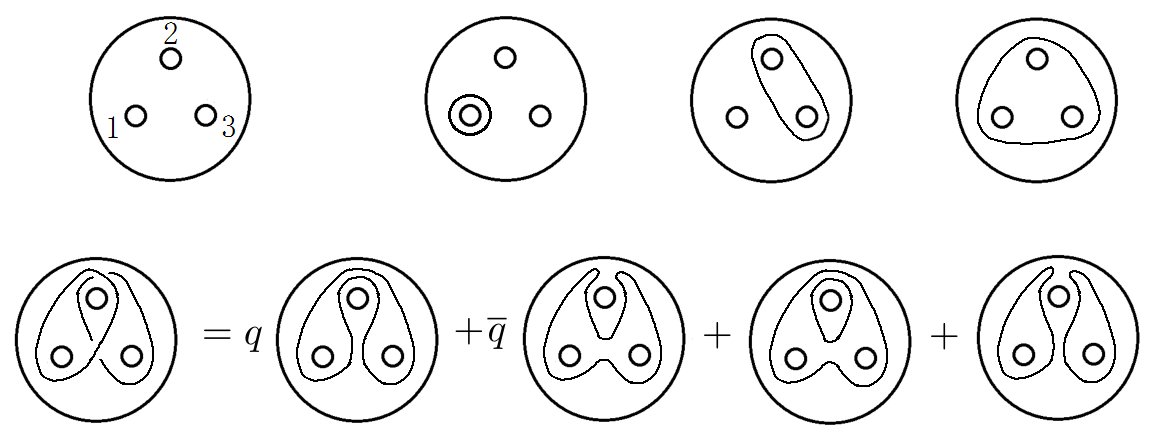}
  \caption{First row (from left to right): $\Sigma_{0,4}$, $t_1$, $t_{23}$, $t_{123}$.
  Second row: using skein relations to expand $t_{12}t_{23}$.}\label{fig:Sigma2}
\end{figure}

Display $\Sigma_{0,4}$ and introduce some elements of $\mathcal{S}_3$ in the first row of Figure \ref{fig:Sigma2}. The meanings for $t_2$, $t_3$, $t_{12}$, $t_{13}$, $t_{123\overline{2}}$ are self-explanatory.

We shall frequently use the following facts: (i) each orientation-preserving embedding $\iota:\Sigma_{0,4}\hookrightarrow\Sigma_{0,5}$ induces a $R$-algebra homomorphism $\iota_\ast:\mathcal{S}_3\to\mathcal{S}_4$; (ii) the map $\tau:\Sigma_{0,n+1}\times[0,1]\to\Sigma_{0,n+1}\times[0,1]$, $(\mathsf{x},c)\mapsto(\mathsf{x},1-c)$ induces an automorphism $\tau_\ast$ of $\mathcal{S}_n$ as an abelian group, satisfying
$\tau_\ast(q^{\frac{1}{2}}x)=\overline{q}^{\frac{1}{2}}\tau_\ast(x)$. See \cite[Proposition 2.2 (1)]{Pr99} and its proof.
Consequently, if $h=0$ is a relation in $\mathcal{S}_3$, then $\iota_\ast(h)=0$ is one in $\mathcal{S}_4$;
if $f=0$ is a relation in $\mathcal{S}_n$, then so is $\tau_\ast(f)=0$, called the {\it mirror image} of $f=0$. Note that $\tau_\ast(f)$ is obtained from $f$ by replacing $q$ with $\overline{q}$ and rewriting each monomial backwards.

The second row of Figure \ref{fig:Sigma2} reads
\begin{align}
t_{12}t_{23}=qt_{123\overline{2}}+\overline{q}t_{13}+t_2t_{123}+t_1t_3.   \label{eq:12-times-23}
\end{align}
Comparing this with its mirror image and eliminating $t_{123\overline{2}}$ yields
\begin{align}
t_{23}t_{12}=\overline{q}^2t_{12}t_{23}+(q-\overline{q}^3)t_{13}+(1-\overline{q}^2)(t_1t_3+t_2t_{123}).  \label{eq:basic-commutator}
\end{align}

\begin{figure}[h]
  \centering
  \includegraphics[width=11.8cm]{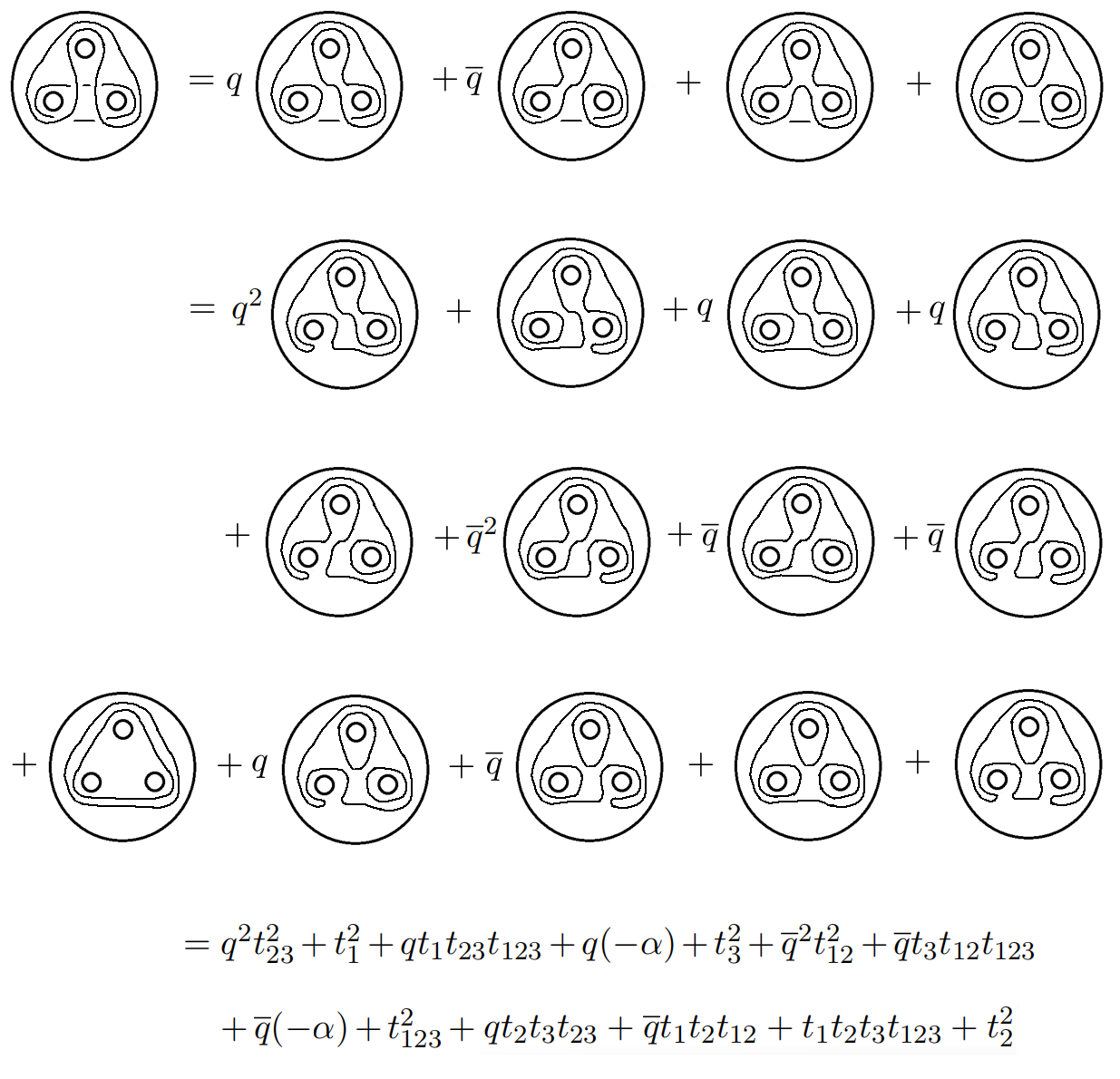}\\
  \caption{Computing $t_{123\overline{2}}t_{13}$ in $\mathcal{S}_3$.}\label{fig:product-1}
\end{figure}

Multiplying (\ref{eq:12-times-23}) by $t_{13}$ on the right, we obtain
\begin{align}
t_{12}t_{23}t_{13}=qt_{123\overline{2}}t_{13}+\overline{q}t_{13}^2+t_2t_{123}t_{13}+t_1t_3t_{13}.   \label{eq:12-times-23-times-13}
\end{align}
As shown in Figure \ref{fig:product-1},
\begin{align*}
t_{123\overline{2}}t_{13}&=t_{123}^2+(qt_1t_{23}+\overline{q}t_3t_{12}+t_1t_2t_3)t_{123}+q^2t_{23}^2+\overline{q}^2t_{12}^2 \nonumber \\
&\ \ \ \ +qt_2t_3t_{23}+\overline{q}t_1t_2t_{12}+t_1^2+t_2^2+t_3^2-\alpha^2.
\end{align*}
Substituting it into (\ref{eq:12-times-23-times-13}) yields
\begin{align}
t_{123}^2&=\overline{q}t_{12}t_{23}t_{13}-(t_1t_2t_3+qt_1t_{23}+\overline{q}t_2t_{13}+\overline{q}t_3t_{12})t_{123}-t_1^2-t_2^2-t_3^2 \nonumber \\
&\ \ \ \ +\alpha^2-qt_2t_3t_{23}-\overline{q}t_1t_3t_{13}-\overline{q}t_1t_2t_{12}-q^2t_{23}^2-\overline{q}^2t_{13}^2-\overline{q}^2t_{12}^2.
\label{eq:t123-square}
\end{align}
This appeared as \cite[Equation (4)]{BP00}, with $A^2=q$, $a_1=t_1$, $a_2=t_3$, $a_3=t_2$, $a_4=t_{123}$, $x_1=t_{13}$, $x_2=t_{23}$, $x_3=t_{12}$.

\begin{figure}[h]
  \centering
  \includegraphics[width=12.5cm]{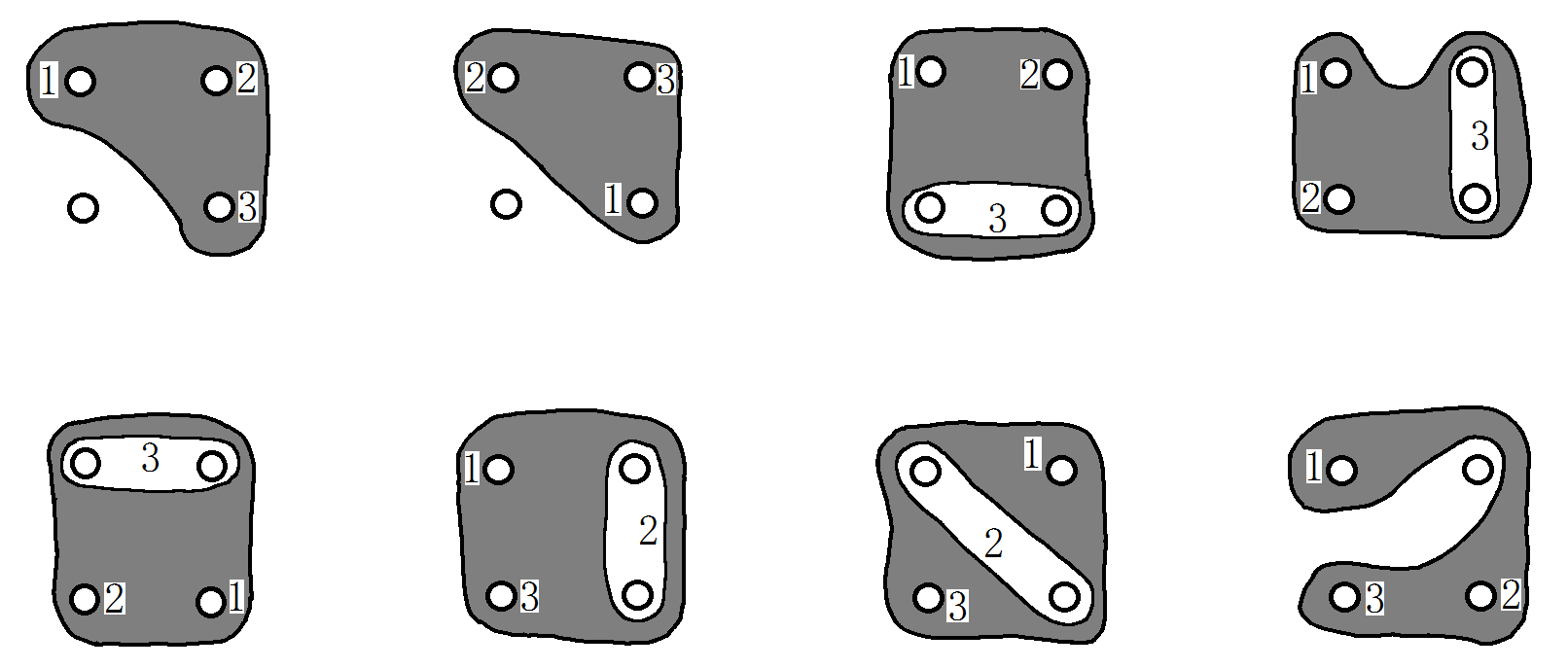}\\
  \caption{First row (from left to right): $\iota_1$--$\iota_4$. Second row: $\iota_5$--$\iota_8$.}\label{fig:embedding}
\end{figure}

To save space, from now on we omit the outer boundary of $\Sigma_{0,5}$.
%Let ${\rm Eq}(t_{13}t_{24})$ denote the relation in $\mathcal{H}$ whose left-hand-side is $t_{13}t_{24}$, and so on.

Let $\iota_1,\ldots,\iota_8:\Sigma_{0,4}\hookrightarrow\Sigma_{0,5}$ denote the orientation-preserving embeddings shown in Figure \ref{fig:embedding}. For each embedding $\iota_i$, we display its image as a shaded region, and label $j$ beside the image of the $j$-th inner boundary circle of $\Sigma_{0,4}$; these determine $\iota_i$ up to isotopy.

We shall use $(\iota_1)_\ast$(\ref{eq:basic-commutator}) to denote the relation in $\mathcal{S}_4$ induced from (\ref{eq:basic-commutator}) by the morphism $(\iota_1)_\ast:\mathcal{S}_3\to\mathcal{S}_4$, and so on. Note that $(\iota_1)_\ast$(\ref{eq:12-times-23}) has the same expression as (\ref{eq:12-times-23}).

As is easy to see, (\ref{eq:commutator-23-12}) $=(\iota_1)_\ast$(\ref{eq:basic-commutator}),
(\ref{eq:commutator-13-12}) $=(\iota_2)_\ast$(\ref{eq:basic-commutator}),
(\ref{eq:reduction-123-123}) $=(\iota_1)_\ast$(\ref{eq:t123-square}), and (\ref{eq:reduction-0-0}) =$(\iota_3)_\ast$(\ref{eq:t123-square}).

\begin{figure}[H]
  \centering
  \includegraphics[width=11cm]{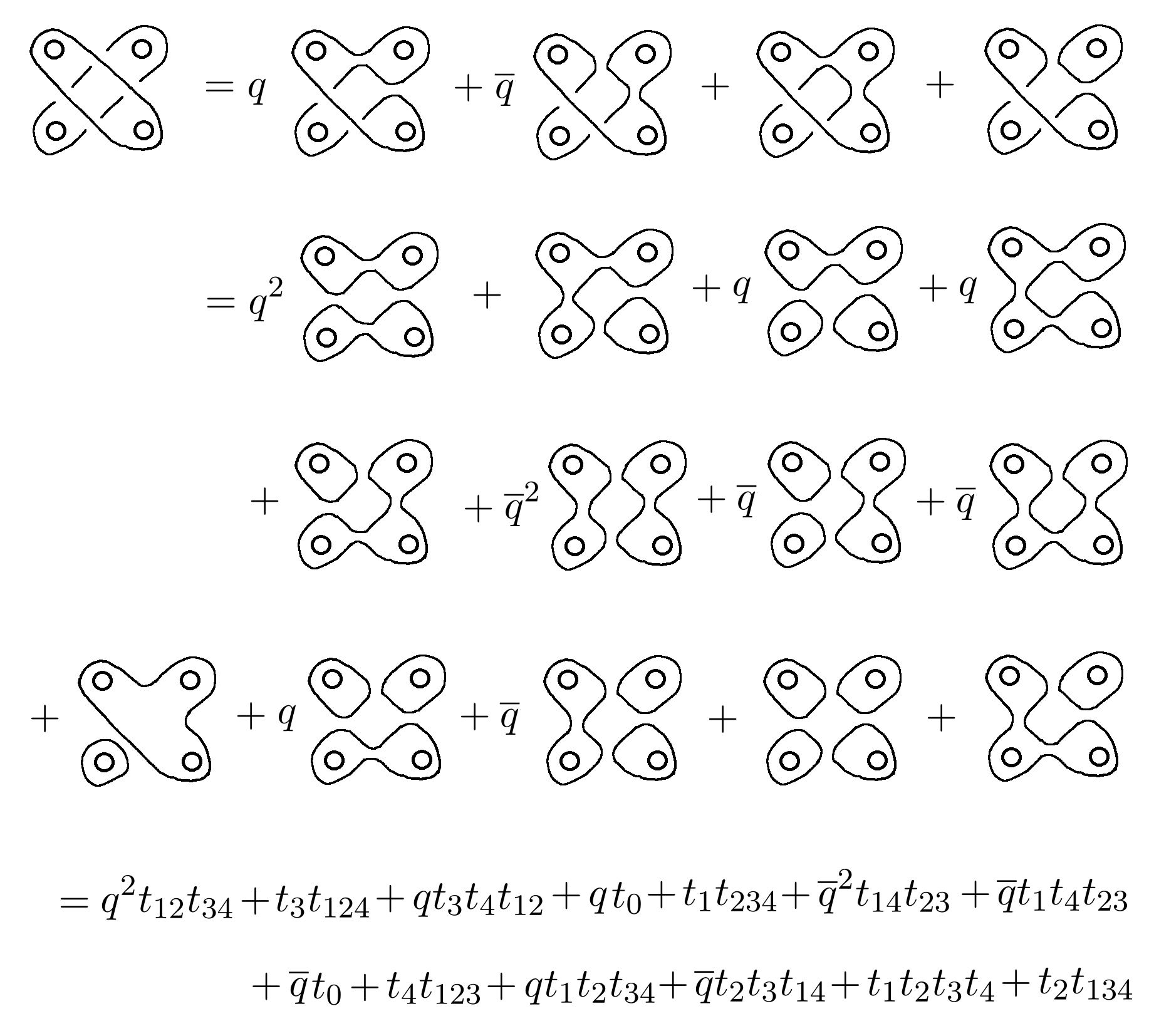}\\
  \caption{Computing $t_{13}t_{24}$.}\label{fig:13-times-24}
\end{figure}

\begin{figure}[h]
  \centering
  \includegraphics[width=11.5cm]{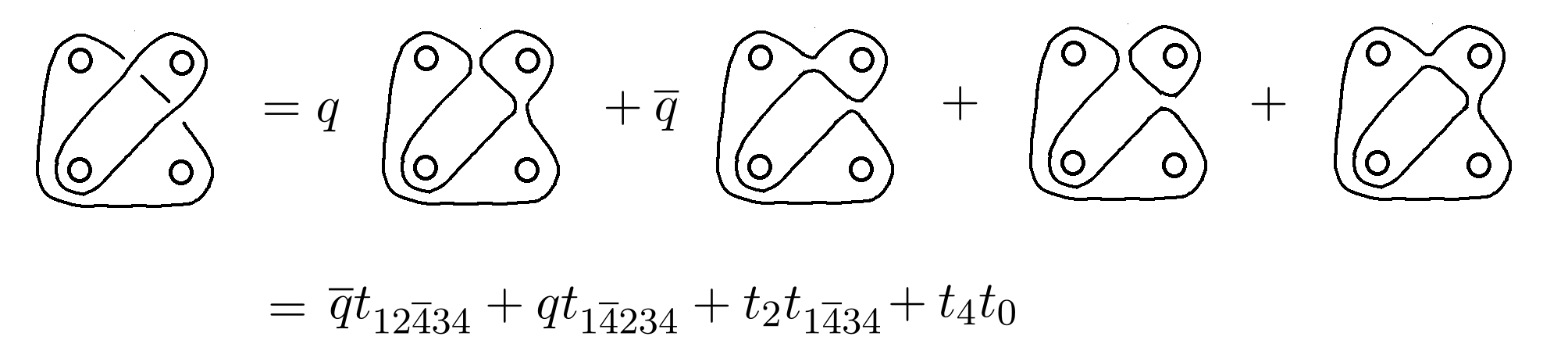}\\
  \caption{Computing $t_{24}t_{134}$. It can also be deduced from (\ref{eq:12-times-23}) using some embedding $\Sigma_{0,4}\hookrightarrow\Sigma_{0,5}$. Here for clarity, we deduce it directly.}\label{fig:product-2}
\end{figure}

\begin{proof}[Proof of (\ref{eq:commutator-24-12})--(\ref{eq:reduction-134-24})]
Similarly as $(\iota_1)_\ast$(\ref{eq:12-times-23}), we can deduce
\begin{align}
t_{12}t_{24}=qt_{124\overline{2}}+\overline{q}t_{14}+t_2t_{124}+t_1t_4,    \label{eq:12-times-24}   \\
t_{34}t_{14}=qt_{1\overline{4}34}+\overline{q}t_{13}+t_4t_{134}+t_1t_3.    \label{eq:34-times-14}
\end{align}
Comparing (\ref{eq:12-times-24}) with its mirror image, we obtain (\ref{eq:commutator-24-12}).

The relations $(\iota_4)_\ast$(\ref{eq:12-times-23}) and $(\iota_5)_\ast$(\ref{eq:12-times-23}) respectively read
\begin{align}
%t_{13}t_{12}&=qt_{12\overline{1}3}+\overline{q}t_{23}+t_1t_{123}+t_2t_3, \label{eq:13-times-12}  \\
t_{14}t_{234}&=qt_{123}+\overline{q}t_{1\overline{4}234}+t_4t_0+t_1t_{23}, \label{eq:14-times-234} \\
t_{34}t_{124}&=qt_{12\overline{4}34}+\overline{q}t_{123}+t_4t_0+t_3t_{12}. \label{eq:34-times-124}
\end{align}
The commutator relation (\ref{eq:commutator-234-14}) results from comparing (\ref{eq:14-times-234}) with its mirror image, and
(\ref{eq:commutator-124-34}) results from comparing (\ref{eq:34-times-124}) with its mirror image.

The reduction relation (\ref{eq:reduction-13-24}) is deduced as in Figure \ref{fig:13-times-24}.

Taking the mirror image of (\ref{eq:34-times-14}), we obtain
\begin{align}
t_{14}t_{34}=\overline{q}t_{1\overline{4}34}+qt_{13}+t_4t_{134}+t_1t_3.   \label{eq:14-times-34}
\end{align}
By Figure \ref{fig:product-2},
\begin{align*}
t_{24}t_{134}=qt_{1\overline{4}234}+\overline{q}t_{12\overline{4}34}+t_2t_{1\overline{4}34}+t_4t_0,   %\label{eq:24-times-134}
\end{align*}
which in combination with (\ref{eq:14-times-234}), (\ref{eq:34-times-124}), (\ref{eq:14-times-34}) yields the reduction relation (\ref{eq:reduction-24-134}).
Its mirror image is (\ref{eq:reduction-134-24}).
\end{proof}

\begin{proof}[Proof of (\ref{eq:reduction-123-234})--(\ref{eq:reduction-0-234})]
As is not difficult to verify, $(\iota_6)_\ast$(\ref{eq:12-times-23}), $(\iota_7)_\ast$(\ref{eq:12-times-23}), $(\iota_8)_\ast$(\ref{eq:12-times-23}) respectively read
\begin{align}
t_{123}t_{234}&=qt_{1234\overline{3}\overline{2}}+\overline{q}t_{14}+t_{23}t_0+t_1t_4,  \label{eq:123-times-234} \\
t_{123}t_{134}&=qt_{234\overline{3}}+\overline{q}t_{12\overline{1}4}+t_{13}t_0+t_2t_4,  \label{eq:123-times-134} \\
t_{123\overline{2}}t_{34}&=qt_{1234\overline{3}\overline{2}}+\overline{q}t_{124\overline{2}}+t_3t_{1234\overline{2}}+t_1t_4,  \label{eq:4-th}
\end{align}

Rotating (\ref{eq:34-times-124}) by $\pi$ yields
\begin{align}
t_{12}t_{234}=qt_{1234\overline{2}}+\overline{q}t_{134}+t_2t_0+t_1t_{34}.   \label{eq:12-times-234}
\end{align}
Multiplying $(\iota_1)_\ast$(\ref{eq:12-times-23}) by $t_{34}$ on the right and using (\ref{eq:4-th}), we obtain
\begin{align*}
t_{12}t_{23}t_{34}&=qt_{123\overline{2}}t_{34}+\overline{q}t_{13}t_{34}+t_2t_{123}t_{34}+t_1t_3t_{34}   \\
&=q\big(qt_{1234\overline{3}\overline{2}}+\overline{q}t_{124\overline{2}}+t_3t_{1234\overline{2}}+t_1t_4\big)
+\overline{q}t_{13}t_{34}+t_2t_{123}t_{34}+t_1t_3t_{34};
\end{align*}
replacing $t_{124\overline{2}}$ and $t_{1234\overline{2}}$ with polynomials respectively via (\ref{eq:12-times-24}) and (\ref{eq:12-times-234}), we are led to
\begin{align}
t_{12}t_{23}t_{34}&=q^2t_{1234\overline{3}\overline{2}}+t_3t_{12}t_{234}+t_2t_{123}t_{34}-\overline{q}t_2t_{124}
-\overline{q}t_3t_{134}+\overline{q}t_{12}t_{24} \nonumber \\
&\ \ \ \ +\overline{q}t_{13}t_{34}-\overline{q}^2t_{14}+(q-\overline{q})t_1t_4-t_2t_3t_0. \label{eq:12-times-23-times-34}
\end{align}
Replacing $t_{1234\overline{3}\overline{2}}$ in (\ref{eq:123-times-234}) via (\ref{eq:12-times-23-times-34}), we obtain
(\ref{eq:reduction-123-234}). Its mirror image is (\ref{eq:reduction-234-123}).

\begin{figure}[h]
  \centering
  \includegraphics[width=11.8cm]{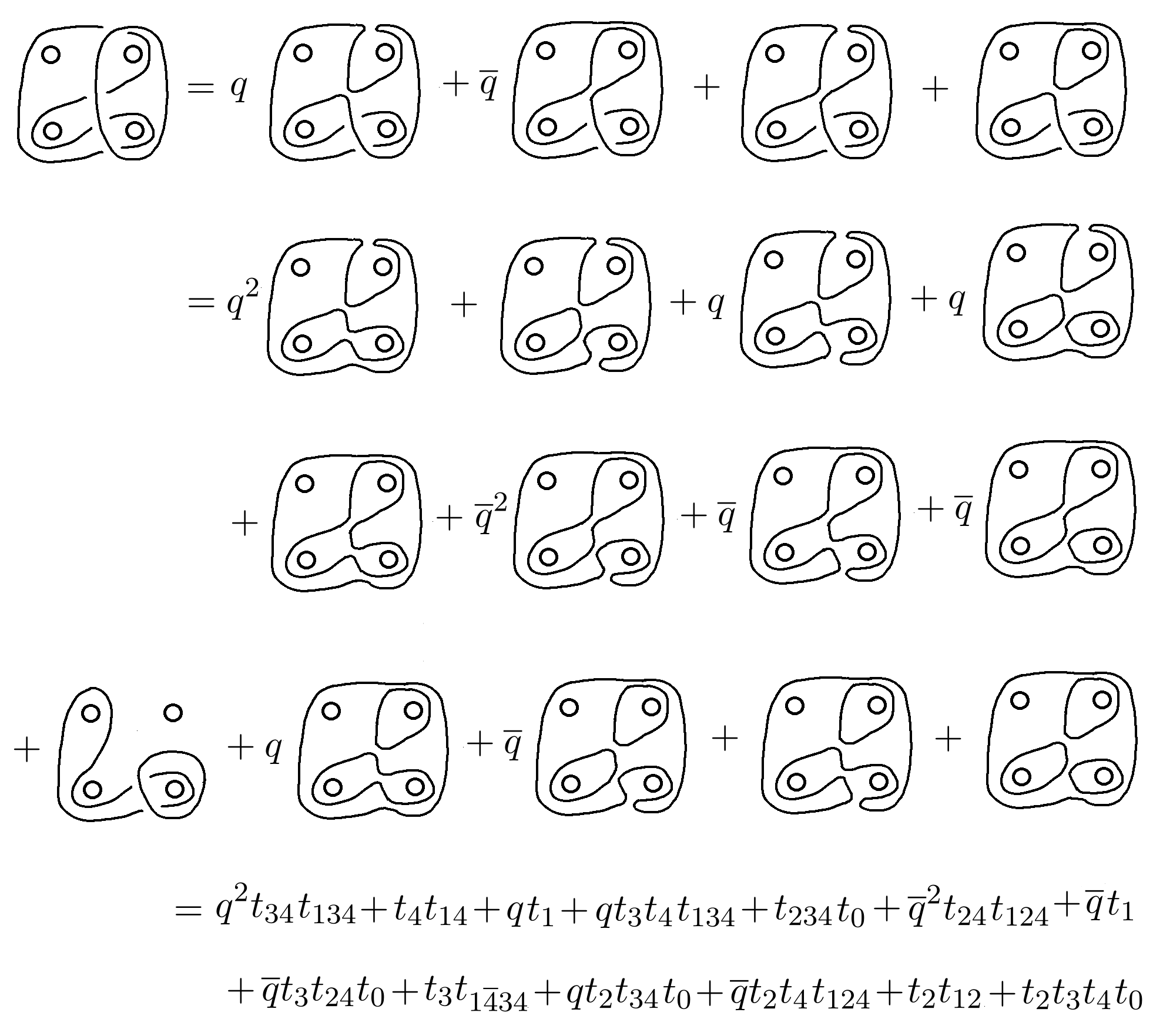}\\
  \caption{Computing $t_{23}t_{12\overline{4}34}$.}\label{fig:product-3}
\end{figure}

Rotating $(\iota_1)_\ast$(\ref{eq:12-times-23}) by $\pi/2$ and taking the mirror image yields
\begin{align*}
t_{12\overline{1}4}=qt_{12}t_{14}-q^2t_{24}-qt_2t_4-qt_1t_{124};
\end{align*}
rotating $(\iota_1)_\ast$(\ref{eq:12-times-23}) by $-\pi/2$ yields
\begin{align*}
t_{234\overline{3}}=\overline{q}t_{23}t_{34}-\overline{q}^2t_{24}-\overline{q}t_2t_4-\overline{q}t_3t_{234}.
\end{align*}
Combining these two equations with (\ref{eq:123-times-134}), we are led to (\ref{eq:reduction-123-134}).

Finally, to verify (\ref{eq:reduction-0-234}), we compute
\begin{align*}
t_{23}t_{12\overline{4}34}=\ &t_3t_{1\overline{4}34}+(t_{234}+qt_2t_{34}+\overline{q}t_3t_{24}+t_2t_3t_4)t_0+(q^2t_{34}+qt_3t_4)t_{134}   \\
&+(\overline{q}^2t_{24}+\overline{q}t_2t_4)t_{124}+t_2t_{12}+t_4t_{14}+\alpha t_1,
\end{align*}
as shown in Figure \ref{fig:product-3}.
Multiplying (\ref{eq:34-times-124}) by $t_{23}$ on the left yields
$$t_{23}t_{34}t_{124}=qt_{23}t_{12\overline{4}34}+\overline{q}t_{24}t_{123}+t_4t_{23}t_0+t_3t_{23}t_{12}.$$
Comparing these equations and using (\ref{eq:14-times-34}) to replace $t_{1\overline{4}34}$ with a polynomial, we obtain (\ref{eq:reduction-0-234}).
\end{proof}

%\begin{rmk}
%\rm Some formulas are indeed complicated. But the key point is the existence of reduction relations and how to obtain them.
%\end{rmk}

\section{On the character variety of the rank $4$ free group} \label{sec:character-variety}

Let $P$ denote the (commutative) polynomial ring over $\mathbb{C}$ generated by $z_i$'s with $1\le i\le n$, $z_{ij}$'s with $1\le i\le j\le n$, and $z_{ijk}$'s with $1\le i<j<k\le n$.

For $i<j$, let $z_{ji}=z_{ij}$. For any $1\le i,j,k\le n$, when $i,j,k$ are pairwise distinct, take a unique permutation $\tau$ with $\tau(i)<\tau(j)<\tau(k)$ and put $z_{ijk}=(-1)^\tau z_{\tau(i)\tau(j)\tau(k)}$; otherwise, put $z_{ijk}=0$.

Let $H_k=\{\underline{a}=(a_1,\ldots,a_k)\colon 1\le a_1<\cdots<a_k\le n\}$ for $k=2,3,4$. For $\underline{a}=(a_1,a_2,a_3)\in H_3$, simplify $z_{a_1a_2a_3}$ as $z_{\underline{a}}$.

For any $\underline{a},\underline{b},\underline{c}\in H_3$ and $1\le b,c,b_i,c_i\le n$, put
\begin{align*}
\theta^{b_1\cdots b_r}_{c_1\cdots c_r}&=\det\big((z_{b_ic_j})_{i,j=1}^r\big),  \\
\xi^{\underline{a}}_{\underline{b}}&=2z_{\underline{a}}z_{\underline{b}}+\theta^{a_1a_2a_3}_{b_1b_2b_3},  \\
\zeta^b_c(\underline{a})&=z_{bc}z_{a_1a_2a_3}-z_{a_1c}z_{ba_2a_3}+z_{a_2c}z_{ba_1a_3}-z_{a_3c}z_{ba_1a_2},  \\
\zeta^{b_1b_2}_{c_1c_2}(\underline{a})&=\theta^{b_1b_2}_{c_1c_2}z_{a_1a_2a_3}-\theta^{b_1a_2}_{c_1c_2}z_{a_1b_2a_3}
+\theta^{b_1a_3}_{c_1c_2}z_{a_1b_2a_2}+\theta^{b_2a_2}_{c_1c_2}z_{a_1b_1a_3}  \\
&\ \ \ -\theta^{b_2a_3}_{c_1c_2}z_{a_1b_1a_2}+\theta^{a_2a_3}_{c_1c_2}z_{a_1b_1b_2},  \\
\zeta^{\underline{b}}_{\underline{c}}(\underline{a})&=\theta^{b_1b_2b_3}_{c_1c_2c_3}z_{a_1a_2a_3}-\theta^{b_1b_2a_3}_{c_1c_2c_3}z_{a_1a_2b_3}
+\theta^{b_1b_3a_3}_{c_1c_2c_3}z_{a_1a_2b_2}-\theta^{b_2b_3a_3}_{c_1c_2c_3}z_{a_1a_2b_1}.
\end{align*}

Let $J$ denote the ideal of $P$ generated by the $\xi^{\underline{a}}_{\underline{b}}$'s for $\underline{a},\underline{b}\in H_3$, and the $\zeta^b_c(\underline{a})$'s for $1\le c\le n$ and $1\le b<a_1<a_2<a_3\le n$.

Call $\xi^{\underline{a}}_{\underline{b}}$ a {\it type I relation}, and call $\zeta^b_c(\underline{a})$ a {\it type II relation}.

The structure of the ring $\mathsf{T}_n$ of ${\rm GL}(2,\mathbb{C})$ invariants of $n$-tuples of $2\times 2$ matrices was known \cite{Dr03}.
Here we state it in a slightly modified form.

Let $\mathsf{M}$ denote the $\mathbb{C}$-algebra of $2\times 2$ matrices. For $\mathbf{x}\in\mathsf{M}$, let
$[[\mathbf{x}]]=\mathbf{x}-(1/2){\rm tr}(\mathbf{x})\cdot\mathbf{e}$, where $\mathbf{e}$ denotes the identity matrix.
Given $1\le i_1,\ldots,i_r\le n$, let $\mathsf{t}_{i_1\cdots i_r}$ denote the function sending
$(\mathbf{x}_1,\ldots,\mathbf{x}_n)\in\mathsf{M}^n$ to ${\rm tr}(\mathbf{x}_{i_1}\cdots\mathbf{x}_{i_r})$, and
let $\mathsf{s}_{i_1\cdots i_r}(\mathbf{x}_1,\ldots,\mathbf{x}_n)={\rm tr}([[\mathbf{x}_{i_1}]]\cdots[[\mathbf{x}_{i_r}]])$.
It is easy to see
\begin{align}
\mathsf{s}_{ij}&=\mathsf{t}_{ij}-\frac{1}{2}\mathsf{t}_i\mathsf{t}_j,  \label{eq:sij}  \\
\mathsf{s}_{ijk}&=\mathsf{t}_{ijk}-\frac{1}{2}(\mathsf{t}_i\mathsf{t}_{jk}+\mathsf{t}_j\mathsf{t}_{ik}+\mathsf{t}_k\mathsf{t}_{ij}
-\mathsf{t}_i\mathsf{t}_j\mathsf{t}_k).  \label{eq:sijk}
\end{align}
Define $\pi:P\to\mathsf{T}_n$ by
$z_i\mapsto\mathsf{t}_i$, $z_{ij}\mapsto\mathsf{s}_{ij}$, $z_{ijk}\mapsto\mathsf{s}_{ijk}.$
By \cite[Theorem 2.3]{Dr03}, $\pi$ is surjective with $\ker(\pi)=J$, so $\mathsf{T}_n\cong P/J$.

Introduce a partial order $\preceq$ on $H=\cup_{k=2}^4H_k$ by
$(a_1,\ldots,a_u)\preceq (b_1,\ldots,b_v)$ if $u\ge v$ and $a_i\le b_i$ for $1\le i\le v$.

Define a variable order by
\begin{align}
z_1<\cdots<z_n&<z_{11}<z_{12}<\cdots<z_{n-1,n}<z_{nn}  \nonumber  \\
&<z_{123}<z_{124}<\cdots<z_{n-2,n-1,n}.  \label{eq:order}
\end{align}

By \cite[Theorem 3.1]{DD01}, with respect to the induced lexicographic monomial order,
${\rm Gr}=\cup_{i=0}^4{\rm Gr}^i$ is a Gr\"obner basis for $J$, where
\begin{align*}
{\rm Gr}^0&=\big\{\theta^{a_1a_2a_3a_4}_{b_1b_2b_3b_4}\colon \underline{a}\preceq\underline{b}\in H_4\big\},  \qquad\qquad
{\rm Gr}^1=\big\{\xi^{\underline{a}}_{\underline{b}}\colon \underline{a},\underline{b}\in H_3\big\}, \\
{\rm Gr}^2&=\big\{\zeta^b_c(\underline{a})\colon  1\le b\le c\le n,\ b<a_1,\ \underline{a}\in H_3\big\}, \\
{\rm Gr}^3&=\big\{\zeta^{b_1b_2}_{c_1c_2}(\underline{a})\colon  \underline{b}\preceq\underline{c}\in H_2, \ \underline{a}\in H_3,\ a_1\le b_1<b_2<a_2\big\},  \\
{\rm Gr}^4&=\big\{\zeta^{\underline{b}}_{\underline{c}}(\underline{a})\colon \underline{a},\underline{b},\underline{c}\in H_3,\ \underline{b}\preceq\underline{c}, \ \underline{a}\preceq(b_1,b_2), \ b_3<a_3\big\}.
\end{align*}
Recall that given $f\in P$, its {\it leading monomial} $L(f)$ is the monomial which is maximal with respect to $<$ among the monomials in $f$. The statement ``${\rm Gr}$ is a Gr\"obner basis for $J$" means that for any $0\ne g\in J$, there exists $h\in {\rm Gr}$ with $L(h)\mid L(g)$.

\begin{rmk}
\rm As a small issue in \cite{DD01}, the $x_{ij}$ with $i>j$ and $u_{i_1i_2i_3}$ without $i_1<i_2<i_3$ were not defined, but they appeared in Theorem 3.1.

This is why we introduce $z_{ij}$ and $z_{ijk}$ for all indices in the second paragraph, according to that $\pi(z_{ij})=\mathsf{s}_{ij}$ satisfies $\mathsf{s}_{ij}=\mathsf{s}_{ji}$, while $\pi(z_{ijk})=\mathsf{s}_{ijk}$ satisfies $\mathsf{s}_{ijk}=\mathsf{s}_{kij}=-\mathsf{s}_{jik}$ (see \cite[Section 2]{Ch24} for instance).
\end{rmk}

Now let $n=4$. Then $H_3=\{(1,2,3),(1,2,4),(1,3,4),(2,3,4)\}$, and it is easy to figure out that
\begin{align*}
{\rm Gr}^0&=\big\{\theta^{1234}_{1234}\big\}, \qquad
{\rm Gr}^1=\big\{\xi^{\underline{a}}_{\underline{b}}\colon \underline{a},\underline{b}\in H_3\big\}, \qquad
{\rm Gr}^2=\big\{\zeta^1_c(234)\colon  1\le c\le 4\big\}, \\
{\rm Gr}^3&=\big\{\zeta^{12}_{c_1c_2}(134)\colon  1\le c_1<c_2\le 4\big\},  \qquad\qquad\ \
{\rm Gr}^4=\big\{\zeta^{123}_{\underline{c}}(124)\colon \underline{c}\in H_3\}.
\end{align*}
The leading monomials of elements in ${\rm Gr}$ are easy to find:
\begin{align*}
L\big(\theta^{1234}_{1234}\big)=z_{11}z_{22}z_{33}z_{44}, \qquad
L\big(\xi^{\underline{a}}_{\underline{b}}\big)=z_{\underline{a}}z_{\underline{b}},  \qquad
L\big(\zeta^1_c(234)\big)=z_{1c}z_{234},  \\
L\big(\zeta^{12}_{c_1c_2}(134)\big)=z_{1c_1}z_{2c_2}z_{134},  \qquad\qquad\qquad
L\big(\zeta^{123}_{\underline{c}}(124)\big)=z_{1c_1}z_{2c_2}z_{3c_3}z_{124}.
\end{align*}
From these it is obvious that
\begin{align}
h\in{\rm Gr}\setminus{\rm Gr}^1\ \Rightarrow\ z_{123}\nmid L(h).  \label{eq:observation}
\end{align}

For $0\ne f\in P$ and a variable $z$, let $\deg_z(f)$ denote the degree of $f$ with respect to $z$.
Let
\begin{align*}
\deg^2(f)&=\deg_{z_{13}}(f)+\deg_{z_{24}}(f),  \\
\deg^3(f)&=\deg_{z_{123}}(f)+\deg_{z_{124}}(f)+\deg_{z_{134}}(f)+\deg_{z_{234}}(f).
\end{align*}
Let $S=\mathbb{C}[z_1,z_2,z_3,z_4,z_{11},z_{22},z_{33},z_{44},z_{12},z_{23},z_{34},z_{14}]$.

\begin{lem}\label{lem:useful}
If $0\ne g\in S[z_{13},z_{24}]\cap J$ with $\deg^3(g)=0$, then $\theta^{1234}_{1234}\mid g$.
\end{lem}

\begin{proof}
Since $g\in J$, by the property of Gr\"obner basis there exists $h\in{\rm Gr}$ such that $L(h)\mid L(g)$.
It is impossible for $h\in\cup_{i=1}^4{\rm Gr}^i$, as otherwise $z_{\underline{a}}\mid L(h)$ for some $\underline{a}\in H_3$, which would imply $z_{\underline{a}}\mid L(g)$. Hence $h=\theta^{1234}_{1234}$.

Let $g_1$ be the remainder when dividing $g$ by $\theta^{1234}_{1234}$, then $L(g_1)<L(h)$ and $\deg^3(g_1)=0$. By the above argument, $L(\theta^{1234}_{1234})\mid L(g_1)$.

Step by step, ultimately we arrive at $\theta^{1234}_{1234}\mid g$.
\end{proof}

\begin{thm} \label{thm:classical-0}
Let $\mathsf{R}=\mathbb{C}[\mathsf{t}_1,\mathsf{t}_2,\mathsf{t}_3,\mathsf{t}_4,\mathsf{s}_{11},\mathsf{s}_{22},\mathsf{s}_{33},\mathsf{s}_{44},
\mathsf{s}_{12},\mathsf{s}_{23},\mathsf{s}_{34},\mathsf{s}_{14}]$.
Then $\mathsf{R}$ is a polynomial ring, and as a $\mathsf{R}$-module, $\mathsf{T}_4$ is freely generated by $\mathsf{E}$, where
\begin{align*}
\mathsf{E}&=\{1,\mathsf{s}_{13},\mathsf{s}_{24},\mathsf{s}_{13}\mathsf{s}_{24}\}
\cup\big\{\mathsf{s}_{13}^k,\mathsf{s}_{13}^k\mathsf{s}_{24},\mathsf{s}_{24}^k,\mathsf{s}_{24}^k\mathsf{s}_{13}\colon k\ge 2\big\}  \\
&\ \ \ \ \cup\big\{\mathsf{s}_{13}^k\mathsf{s}_{123},\mathsf{s}_{13}^k\mathsf{s}_{134},\mathsf{s}_{24}^k\mathsf{s}_{124},
\mathsf{s}_{24}^k\mathsf{s}_{234}\colon k\ge 0\big\}.
\end{align*}
\end{thm}

\begin{proof}
{\bf Step 1}. To show $\mathsf{R}$ is a polynomial ring, it suffices to show $S\cap J=0$.

Assume there exists $0\ne g\in S\cap J$. Then $\deg^2(g)=\deg^3(g)=0$. By Lemma \ref{lem:useful}, $\theta^{1234}_{1234}\mid g$.
However, $\deg^2(\theta^{1234}_{1234})=4$. This is a contradiction.

\medskip

{\bf Step 2}. Suppose $m$ is an arbitrary monomial in $z_{13}$, $z_{24}$, $z_{123}$, $z_{124}$, $z_{134}$, $z_{234}$.
We show $m\in S\langle E\rangle+J$, where
\begin{align*}
E&=\{1,z_{13},z_{24},z_{13}z_{24}\}
\cup\big\{z_{13}^k,z_{13}^kz_{24},z_{24}^k,z_{24}^kz_{13}\colon k\ge 2\big\}  \\
&\ \ \ \ \cup\big\{z_{13}^kz_{123},z_{13}^kz_{134},z_{24}^kz_{124},z_{24}^kz_{234}\colon k\ge 0\big\}.
\end{align*}

If $\deg^3(m)\ge 2$, then modulo type I relations, $m$ can be replaced by a $S$-linear combination of some monomials $m_i$ with $\deg^3(m_i)=\deg^3(m)-2$. Hence we only need to consider the case $\deg^3(m)\le 1$.

If $\deg^3(m)=0$, i.e., $m=z_{13}^{k_1}z_{24}^{k_2}$ for some $k_1,k_2\ge 0$, then the remainder of dividing $m$ by $\theta^{1234}_{1234}$ belongs to $S\langle E\rangle$. Hence $m\in S\langle E\rangle+J$.

Now assume $\deg^3(m)=1$. We use induction on $\deg^2(m)$ to prove $m\in S\langle E\rangle+J$, which obviously holds when $\deg^2(m)=0$. Suppose $\deg^2(m)>0$ and that $m'\in S\langle E\rangle+J$ for each $m'$ with $\deg^2(m')<\deg^2(m)$.

Assume $\deg_{z_{123}}(m)=1$, so that $\deg_{x}(m)=0$ for all $x\in\{z_{124},z_{134},z_{234}\}$.
If $z_{24}\nmid m$, then $m=z_{13}^kz_{123}$ for some $k\ge 0$, so $m\in S\langle E\rangle$.
Otherwise, modulo $\zeta^1_2(234)$ we can replace $z_{24}z_{123}$ with $z_{12}z_{234}-z_{22}z_{134}+z_{23}z_{124}$, so as to replace $m$ with $m_1-m_2+m_3$ for some $m_i$ with $\deg^2(m_i)<\deg^2(m)$.
By the inductive hypothesis, $m_i\in S\langle E\rangle+J$. Hence $m\in S\langle E\rangle+J$.

Similarly, $m\in S\langle E\rangle+J$ when $\deg_{x}(m)=1$ for $x\in\{z_{134},z_{124},z_{234}\}$.

Consequently, each monomial in $\mathsf{s}_{13}$, $\mathsf{s}_{24}$, $\mathsf{s}_{123}$, $\mathsf{s}_{124}$, $\mathsf{s}_{134}$, $\mathsf{s}_{234}$ belongs to $\mathsf{R}\langle\mathsf{E}\rangle$. Thus, $\mathsf{T}_4=\mathsf{R}\langle\mathsf{E}\rangle$.

\medskip

{\bf Step 3}. Assume there is a $\mathsf{R}$-linearly dependence among the elements of $\mathsf{E}$.
Then there exist $g_0\in S[z_{13},z_{24}]$, $g_1,g_3\in S[z_{24}]$, $g_2,g_4\in S[z_{13}]$ such that
\begin{align}
g_0+g_1z_{234}+g_2z_{134}+g_3z_{124}+g_4z_{123}\in J.  \label{eq:dependent}
\end{align}
Dividing by $\theta^{1234}_{1234}$ if necessary, we can assume
\begin{align}
\text{no\ monomial\ in\ }g_0\ \text{is\ divisible\ by\ }z_{13}^2z_{24}^2.  \label{eq:condition}
\end{align}

Multiplying (\ref{eq:dependent}) by $z_{123}$ and using type I relations to replace $z_{ijk}z_{123}$ with $(-1/2)\theta^{ijk}_{123}$,
we obtain $w:=g_0z_{123}+v\in J$, with
$$v=-\frac{1}{2}\big(g_1\theta^{234}_{123}+g_2\theta^{134}_{123}+g_3\theta^{124}_{123}+g_4\theta^{123}_{123}\big).$$
By the property of Gr\"obner basis, $L(h)\mid L(w)$ for some $h\in{\rm Gr}$.
Recalling (\ref{eq:order}), due to the presence of $z_{123}$, we have $m<L(g_0z_{123})$ for each monomial $m$ in $v$, so that $L(w)=L(g_0z_{123})$.
Obviously,
$L(\xi^{\underline{a}}_{\underline{b}})\nmid L(g_0z_{123})$ for any $\underline{a},\underline{b}$, so $h\notin{\rm Gr}^1$.
By (\ref{eq:observation}), actually $L(h)\mid L(g_0)$, which forces $h=\theta^{1234}_{1234}$.

Let $r$ denote the remainder of dividing $g_0$ by $h$, then $L(r)<L(g_0)$ and $rz_{123}+v\in J$. By the above argument, $L(\theta^{1234}_{1234})\mid L(r)$.

Step by step, eventually we arrive at $\theta^{1234}_{1234}\mid g_0$. By (\ref{eq:condition}), $g_0=0$, so $v\in J$.
By Lemma \ref{lem:useful}, $\theta^{1234}_{1234}\mid v$.

Recall that $g_1,g_3\in S[z_{24}]$ and $g_2,g_4\in S[z_{13}]$.
As is easy to see,
\begin{align*}
\deg_{z_{13}}(\theta^{234}_{123})=\deg_{z_{13}}(\theta^{124}_{123})=1, \qquad
\deg_{z_{24}}(\theta^{134}_{123})=1, \quad \deg_{z_{24}}(\theta^{123}_{123})=0.
\end{align*}
Hence no monomial in $v$ is divisible by $z_{13}^2z_{24}^2$.
This forces $v=0$.

It then follows that $g_2=g_4=0$. Indeed, note that
$$\deg_{z_{13}}(\theta^{134}_{123})=\deg_{z_{13}}(\theta^{123}_{123})=2.$$
For $i\in\{2,4\}$, when $g_i\ne 0$, suppose the highest order term of $g_i$ with respect to $z_{13}$ is $e_iz_{13}^{n_i}$, with $0\ne e_i\in S$ and $n_i\ge0$; when $g_i=0$, put $e_i=0$, $n_i=0$.
If $g_2\ne 0$ or $g_4\ne 0$, then the highest order term of $v$ with respect to $z_{13}$ would be
$$\begin{cases}
e_2z_{24}z_{13}^{n_2+2}, &n_2>n_4 \\
(e_2z_{24}+e_4z_{22})z_{13}^{n_2+2},&n_2=n_4 \\
e_4z_{22}z_{13}^{n_4+2},&n_2<n_4
\end{cases},$$
which is necessarily nonzero. Thus, $g_2=g_4=0$.

Multiplying (\ref{eq:dependent}) by $z_{234}$ and using type I relations to replace $z_{ijk}z_{234}$ with $(-1/2)\theta^{ijk}_{234}$, we obtain
$u:=g_1\theta^{234}_{234}+g_3\theta^{124}_{234}\in J.$
By Lemma \ref{lem:useful}, $\theta^{1234}_{1234}\mid u$.
However, as is easy to see,
$$\deg_{z_{13}}(\theta^{234}_{234})=0, \qquad   \deg_{z_{13}}(\theta^{124}_{234})=1.$$
Hence no monomial in $u$ is divisible by $z_{13}^2z_{24}^2$.
This forces $u=0$.

It then follows that $g_1=g_3=0$. Indeed, note that
$$\deg_{z_{24}}(\theta^{234}_{234})=\deg_{z_{24}}(\theta^{124}_{234})=2.$$
For $i=1,3$, when $g_i\ne 0$, suppose the highest order term of $g_i$ with respect to $z_{24}$ is $e_iz_{24}^{n_i}$, with $0\ne e_i\in S$ and $n_i\ge0$; when $g_i=0$, put $e_i=0$, $n_i=0$.
If $g_1\ne 0$ or $g_3\ne 0$, then the highest order term of $u$ with respect to $z_{24}$ would be
$$\begin{cases}
e_1z_{33}z_{24}^{n_1+2}, &n_1>n_3 \\
(e_1z_{33}+e_3z_{13})z_{24}^{n_1+2},&n_1=n_3 \\
e_3z_{13}z_{24}^{n_3+2},&n_1<n_3
\end{cases},$$
which is necessarily nonzero. Thus, $g_1=g_3=0$.

Therefore, the elements of $\mathsf{E}$ are $\mathsf{R}$-linearly independent.
\end{proof}

Let $\Gamma$ be a finitely presented group. 
Given $\rho\in\hom(\Gamma,{\rm SL}(2,\mathbb{C}))$, its {\it character} is the function $\chi_\rho:\Gamma\to\mathbb{C}$, $x\mapsto{\rm tr}(\rho(x))$.
The ${\rm SL}(2,\mathbb{C})$-{\it character variety} $\mathcal{X}(\Gamma)$ is defined as the GIT quotient
of $\hom(\Gamma,{\rm SL}(2,\mathbb{C}))$ by ${\rm SL}(2,\mathbb{C})$ which acts via conjugation. As an affine algebraic set,
$$\mathcal{X}(\Gamma)=\{\chi_\rho\colon\rho\in\hom(\Gamma,{\rm SL}(2,\mathbb{C}))\}.$$

In the case $\Gamma=F_4$, the free group of rank $4$, it is known (see \cite[Section 2]{ABL18}) that the coordinate ring $\mathbb{C}[\mathcal{X}(F_4)]$ is isomorphic to the quotient of $\mathsf{T}_4$ by the ideal generated by $\mathsf{s}_{ii}-(1/2)\mathsf{t}_i^2+2$ for $1\le i\le 4$. Actually, the relation $\mathsf{s}_{ii}=(1/2)\mathsf{t}_i^2-2$ comes from imposing the condition $\det(\mathbf{x}_i)=1$.

Observe that $\mathbb{C}[\mathcal{X}(F_4)]$ is also a subring of $\mathsf{T}_4$; in other words, $\mathbb{C}[\mathcal{X}(F_4)]$ is a direct summand of $\mathsf{T}_4$.
Recall (\ref{eq:sij}), (\ref{eq:sijk}). The following is an immediate corollary of Theorem \ref{thm:classical-0}:
\begin{thm}\label{thm:classical}
Let $\mathsf{Q}=\mathbb{C}[\mathsf{t}_1,\mathsf{t}_2,\mathsf{t}_3,\mathsf{t}_4,\mathsf{t}_{12},\mathsf{t}_{23},\mathsf{t}_{34},\mathsf{t}_{14}]$.
Then $\mathsf{Q}$ is a polynomial ring, and $\mathbb{C}[\mathcal{X}(F_4)]$ is a free $\mathsf{Q}$-module generated by $\mathsf{A}$, where
\begin{align*}
\mathsf{A}&=\{1,\mathsf{t}_{13},\mathsf{t}_{24},\mathsf{t}_{13}\mathsf{t}_{24}\}
\cup\big\{\mathsf{t}_{13}^k,\mathsf{t}_{13}^k\mathsf{t}_{24},\mathsf{t}_{24}^k,\mathsf{t}_{24}^k\mathsf{t}_{13}\colon k\ge 2\big\}  \\
&\ \ \ \ \cup\big\{\mathsf{t}_{13}^k\mathsf{t}_{123},\mathsf{t}_{13}^k\mathsf{t}_{134},\mathsf{t}_{24}^k\mathsf{t}_{124},
\mathsf{t}_{24}^k\mathsf{t}_{234}\colon k\ge 0\big\}.
\end{align*}
Consequently, $\mathsf{B}$ is a basis for $\mathbb{C}[\mathcal{X}(F_4)]$ as a vector space over $\mathbb{C}$, where
$$\mathsf{B}=\big\{\mathsf{t}_1^{i_1}\mathsf{t}_2^{i_2}\mathsf{t}_3^{i_3}\mathsf{t}_4^{i_4}\mathsf{t}_{12}^{i_5}\mathsf{t}_{23}^{i_6}
\mathsf{t}_{34}^{i_7}\mathsf{t}_{14}^{i_8}\mathsf{a}\colon i_1,\ldots,i_8\ge 0, \ \mathsf{a}\in\mathsf{A}\big\}$$
\end{thm}

\section{A monomial basis for $\mathcal{S}_4$}

In this section, let $R=\mathbb{Z}[q^{\pm\frac{1}{2}}]$.

For the commutative ring $\mathbb{C}$, take $-1$ as the fixed invertible element.
For each $3$-manifold $M$, by \cite{Bu97} there is a surjective ring homomorphism
\begin{align*}
\phi:\mathcal{S}(M;\mathbb{C})\twoheadrightarrow\mathbb{C}[\mathcal{X}(\pi_1(M))].  %\label{eq:classicalization}
\end{align*}
Explicitly, $\phi$ takes a link $\sqcup_{i=1}^m\mathbf{k}_i\hookrightarrow M$ (the $\mathbf{k}_i$'s being components) to the function
$\mathcal{X}(\pi_1(M))\to\mathbb{C}$ sending a character $\chi$ to $\prod_{i=1}^m(-\chi([\mathbf{k}_i]))$,
where $[\mathbf{k}_i]$ denotes the conjugacy class determined by $\mathbf{k}_i$.

Let $\epsilon: \mathcal{S}(M):=\mathcal{S}(M;R)\to\mathbb{C}[\mathcal{X}(\pi_1(M))]$ denote the composite of $\phi$ and the map
$\mathcal{S}(M)\to\mathcal{S}(M;\mathbb{C})$ induced by setting $q^{\frac{1}{2}}$ to $-1$.

\begin{lem}\label{lem:strategy}
Suppose $\mathcal{S}(M)$ is torsion-free as a module over $R$. If $\mathcal{B}$ is a generating set for $\mathcal{S}(M)$ such that $\epsilon(\mathcal{B})$ is $\mathbb{C}$-linearly independent, then $\mathcal{B}$ is a basis for $\mathcal{S}(M)$.
\end{lem}

\begin{proof}
Assume there is a nontrivial relation ${\sum}_i\mu_ib_i=0$, with $\mu_i\in R$, $b_i\in\mathcal{B}$.

Let $k_0$ be the maximal $k$ such that $(q^{\frac{1}{2}}+1)^k\mid\mu_i$ for all $i$; write $\mu_i$ as $(q^{\frac{1}{2}}+1)^{k_0}\mu'_i$.
Let $\lambda_i$ be the evaluation of $\mu'_i$ at $-1$. Then $\lambda_{i}\ne 0$ for some $i$.

Since $\mathcal{S}(M)$ is torsion-free, we have ${\sum}_i\mu'_ib_i=0$. Hence ${\sum}_i\lambda_i\epsilon(b_i)=0$. This contradicts the assumption that $\epsilon(\mathcal{B})$ is $\mathbb{C}$-linearly independent.

Thus, $\mathcal{B}$ is $R$-linearly independent, i.e. $\mathcal{B}$ is a basis for $\mathcal{S}(M)$.
\end{proof}

Recall $\mathcal{G}=\{t_{i_1\cdots i_r}\colon 1\le i_1<\cdots<i_r\le 4\}$, and $t_0=t_{1234}$.
Call $t_1$, $t_2$, $t_3$, $t_4$, $t_{12}$, $t_{23}$, $t_{34}$, $t_{14}$ {\it light}, and call $t_0$, $t_{13}$, $t_{24}$, $t_{123}$, $t_{124}$, $t_{134}$, $t_{234}$ {\it heavy}.

Set $|t_{i_1\cdots i_r}|=r$; in particular, $|t_0|=4$. Set $|1|=0$ as a convention.
Given a monomial $u=g_1\cdots g_m$ with $g_i\in\mathcal{G}$, define $|u|={\sum}_{i=1}^m|g_i|$,
%$$|u|={\sum}_{i=1}^m|g_i|, \qquad  |u|_1=|u'|, \qquad   |u|_2=|u''|,$$
and put
$$\|u\|=(|u|,|\dot{u}|,|\ddot{u}|),$$
where $\dot{u}$ is the product of the $g_i$'s with $|g_i|\ge 2$, and $\ddot{u}$ is the product of heavy generators in $u$.

Given monomials $u,v$, write $\|v\|\prec\|u\|$ and say that $v$ is {\it simpler} than $u$ if 
one of the following cases holds: (i) $|v|<|u|$; (ii) $|v|=|u|$, $|\dot{v}|<|\dot{u}|$; (iii) $|v|=|u|$, $|\dot{v}|=|\dot{u}|$, $|\ddot{v}|<|\ddot{u}|$. For instance, $t_1t_{234}$, $t_{12}t_{34}$ are simpler than $t_{13}t_{24}$;
$\|t_{234}t_{14}\|\prec\|t_{24}t_{134}\|$, and $\|t_{12}t_{234}t_{134}\|\prec\|t_0^2\|$.

For $\beta\in R$ and monomials $u,u'$, write $u\equiv \beta u'$ and say that $u$ can be {\it displaced} by $\beta u'$, if $\|u\|=\|u'\|$ and $u-\beta u'$ equals a $R$-linear combination of monomials simpler than $u$. Write $u\sim 0$ and say that $u$ can be {\it reduced}, if $u$ equals a $R$-linear combination of monomials simpler than $u$.

\begin{lem}\label{lem:reduce}
Each of $t_0^2$, $t_0t_{234}$, $t_{24}t_{134}$, $t_{134}t_{24}$, $t_{134}^2$, $t_{134}t_{124}$, $t_{134}t_{234}$ can be reduced;
$t_{13}t_{24}\equiv t_{24}t_{13}\equiv \alpha t_0$, and $t_{134}t_{123}\equiv t_{13}t_0$.
\end{lem}

\begin{proof}
Looking into the reduction relations carefully, we have the following:
\begin{enumerate}
  \item By (\ref{eq:reduction-0-0}), $t_0^2\sim 0$; by (\ref{eq:reduction-0-234}), $t_0t_{234}\sim 0$.
  \item By (\ref{eq:reduction-24-134}), $t_{24}t_{134}\sim 0$; by (\ref{eq:reduction-134-24}), $t_{134}t_{24}\sim 0$.
  \item By (\ref{eq:reduction-123-123}), $t_{123}^2\sim 0$; rotating by $\pi$ yields $t_{134}^2\sim 0$.
  \item By (\ref{eq:reduction-123-234}), $t_{123}t_{234}\sim 0$; rotating it by $\pi$ yields $t_{134}t_{124}\sim 0$.
  \item By (\ref{eq:reduction-234-123}), $t_{234}t_{123}\sim 0$; rotating it by $-\pi/2$ yields $t_{134}t_{234}\sim 0$.
  \item By (\ref{eq:reduction-13-24}), $t_{13}t_{24}\equiv \alpha t_0$; rotating it by $\pi/2$ yields
        $t_{24}t_{13}\equiv \alpha t_0$.
  \item By (\ref{eq:reduction-123-134}), $t_{123}t_{134}\equiv t_{13}t_0$; rotating it by $\pi$ yields $t_{134}t_{123}\equiv t_{13}t_0$.
\end{enumerate}
\end{proof}

\begin{proof}[Proof of Theorem \ref{thm:main-1}]
The assertion is equivalent to that
$$\mathcal{C}=\big\{t_1^{i_1}t_2^{i_2}t_3^{i_3}t_4^{i_4}t_{12}^{i_5}t_{23}^{i_6}t_{34}^{i_7}t_{14}^{i_8}a
\colon i_1,\ldots,i_8\ge 0,\ a\in\mathcal{A}\big\}$$
is a basis for $\mathcal{S}_4$ as a $R$-module; recall that
\begin{align*}
\mathcal{A}=\ &\{1,t_0,t_{123},t_{124},t_{134},t_{234}\}  \\
&\cup\{t_{13}^k,\  t_{13}^kt_0, \ t_{13}^kt_{123},\ t_{13}^kt_{134}, \ t_{24}^k, \ t_{24}^kt_0, \ t_{24}^kt_{124}, \ t_{24}^kt_{234} \colon k\ge 1\}.
\end{align*}

By \cite[Theorem 2.3 (b)]{Pr99}, $\mathcal{S}_4$ is free, so it is torsion-free.

By definition, $\epsilon(t_{i_1\cdots i_r})=-\mathsf{t}_{i_1\cdots i_r}$.
In virtue of (\ref{eq:reduction-13-24}), it is not difficult to see that $\epsilon(\mathcal{C})$ is related to $\mathsf{B}$ via an invertible linear map.
By Lemma \ref{lem:strategy} and Theorem \ref{thm:classical}, it suffices to show $\mathcal{S}_4=R\langle\mathcal{C}\rangle$.

Observe that modulo simpler terms and up to $q^{\pm 2}$, we can interchange a light generator $g_0$ with any $g\in\mathcal{G}$.
This is obvious when $g_0\in\{t_1,\ldots,t_4\}$. When $g_0=t_{12}$, either $gt_{12}=t_{12}g$ or $g\in\{t_{23},t_{14},t_{13},t_{24},t_{134},t_{234}\}$; in the latter case, from the commutator relations we can see $gt_{12}\equiv q^{2\varepsilon}t_{12}g$ for some $\varepsilon\in\{\pm1\}$. The situation is similar when
$g_0\in\{t_{23},t_{34},t_{14}\}$.

%In other words, modulo simpler terms and up to an invertible factor, we can put $g_0$ at any prescribed place.
%Consequently, $g_0\mathcal{C}\subset R\langle\mathcal{C}\rangle$.

We prove $u\in R\langle\mathcal{C}\rangle$ for monomials $u$ by induction on $\|u\|$. Clearly, $u\in R\langle\mathcal{C}\rangle$ if $|u|\le 3$.
Suppose $|u|>3$ and that $u'\in R\langle\mathcal{C}\rangle$ for all monomials $u'$ simpler than $u$.
Let $u=g_1\cdots g_m$ with $g_i\in\mathcal{G}$.

If some $g_i$ is light, then by the above observation, we can move $g_i$ to the leftmost place, so
$u\equiv q^{2\mu}g_iv$ for some $\mu\in\mathbb{Z}$, with $v=g_1\cdots\widehat{g_i}\cdots g_m$ (omitting $g_i$).
By the inductive hypothesis, $u-q^{2\mu}g_iv\in R\langle\mathcal{C}\rangle$ and $v\in R\langle\mathcal{C}\rangle$.
Considering a general element of $\mathcal{C}$, for each $a\in\mathcal{A}$, we have
$$g_it_1^{i_1}t_2^{i_2}t_3^{i_3}t_4^{i_4}t_{12}^{i_5}t_{23}^{i_6}t_{34}^{i_7}t_{14}^{i_8}a
\equiv q^{2\nu}t_1^{i_1}t_2^{i_2}t_3^{i_3}t_4^{i_4}t_{12}^{i'_5}t_{23}^{i'_6}t_{34}^{i'_7}t_{14}^{i'_8}a,$$
for some $\nu\in\mathbb{Z}$, where $i'_5=i_5+1$ if $g_i=t_{12}$ and $i'_5=i_5$ otherwise, and so forth.
Hence by the inductive hypothesis again, $g_iv\in R\langle\mathcal{C}\rangle$. Thus, $u\in R\langle\mathcal{C}\rangle$.

%Consequently, each monomial can be written as a $R$-linear combination of monomials of the form $t_1^{i_1}t_2^{i_2}t_3^{i_3}t_4^{i_4}t_{12}^{i_5}t_{23}^{i_6}t_{34}^{i_7}t_{14}^{i_8}w$, where $w$ is a product of heavy generators.
%It suffices to show $w\in R\langle\mathcal{C}\rangle$.

Now suppose each $g_i$ is heavy. We shall apply Lemma \ref{lem:reduce} repeatedly.
\begin{enumerate}
  \item Whenever $g_i=g_j=t_0$ for some $i<j$, since $t_0$ is central, we can move $g_i$ to meet $g_j$ and reduce $t_0^2$, so as to reduce $u$.
  \item Suppose $g_i=t_0$ for exactly one $i$. Just assume $g_1=t_0$.
  \begin{enumerate}
    \item When $g_j=t_{234}$ for some $j$, we can move $t_0$ to meet $t_{234}$, and then reduce $t_0t_{234}$, so as to reduce $u$.

          Similarly for the cases when some $g_j\in\{t_{123},t_{124},t_{134}\}$.
    \item Otherwise, $|g_j|=2$ for all $j>1$. %, i.e. $g_j\in\{t_{13},t_{24}\}$.
          If $u=t_0t_{13}^{m-1}$ or $u=t_0t_{24}^{m-1}$, we are done;
          otherwise, $\{g_j,g_{j+1}\}=\{t_{13},t_{24}\}$ for some $j$, we can displace $g_jg_{j+1}$ by $\alpha t_0$, and then reduce $t_0^2$ as in Case 1.
  \end{enumerate}
  \item Suppose $\max\{|g_1|,\ldots,|g_m|\}=3$. Let $i_0=\min\{i\colon |g_i|=3\}$. Just assume $g_{i_0}=t_{134}$;
        the cases $g_{i_0}\in\{t_{123},t_{124},t_{234}\}$ are similar.
  \begin{enumerate}
    \item If there exists $j<i_0$ such that $g_j=t_{24}$ and $g_k=t_{13}$ for all $j<k<i_0$, then
          $g_j\cdots g_{i_0}=t_{24}t_{13}^{i_0-j-1}t_{134}=t_{24}t_{134}t_{13}^{i_0-j-1}$, which can be reduced, so can $u$.
    \item If $g_j=t_{13}$ for all $j\ne i_0$, then $u=t_{13}^{m-1}t_{134}$, and we are done.
    \item Otherwise, $g_j=t_{13}$ for all $j<i_0$ and $g_k\ne t_{13}$ for some $k>i_0$. Just assume $g_{i_0+1}\ne t_{13}$; otherwise we can
          simply interchange $g_{i_0}$ with $g_{i_0+1}$ and go on, with $i_0+1$ in place of $i_0$.

          When $g_{i_0+1}=t_{123}$, we can displace $g_{i_0}g_{i_0+1}$ by $t_{13}t_0$, and turn to Case 2;
          when $g_{i_0+1}\ne t_{123}$, i.e. $g_{i_0+1}\in\{t_{24},t_{124},t_{134},t_{234}\}$, we can reduce $g_{i_0}g_{i_0+1}$.
  \end{enumerate}
  \item Suppose all $g_i\in\{t_{13},t_{24}\}$. If $\{g_j,g_{j+1}\}=\{t_{13},t_{24}\}$ for some $j$, then %as in Case 2 (b),
        we can displace $g_jg_{j+1}$ by $\alpha t_0$, and turn to Case 2.
        Otherwise, $u=t_{13}^{m}$ or $u=t_{24}^{m}$.
\end{enumerate}

The proof is completed.
\end{proof}

As noticed in Section 1, only the relations in $\mathcal{H}$ are used to transform each monomial into $R\langle\mathcal{C}\rangle$, so Theorem \ref{thm:main-2} has been proved.

\bigskip

\noindent
{\bf Acknowledgement}

The author would like to deeply thank the referees for carefully reading the manuscript and giving many valuable suggestions and comments, so that the paper can be largely improved.

\medskip

\noindent
Haimiao Chen (orcid: 0000-0001-8194-1264)\ \ \  \emph{chenhm@math.pku.edu.cn} \\
Department of Mathematics, Beijing Technology and Business University, \\
Liangxiang Higher Education Park, Fangshan District, Beijing, China.

\end{document}